\documentclass[12pt]{article}
\usepackage[utf8]{inputenc}

\usepackage{amsmath}
\usepackage{amsfonts}
\usepackage{amssymb}
\usepackage{graphicx}
\usepackage{indentfirst}
\usepackage{mathdots}
\usepackage{amsmath}
\usepackage{amssymb}
\usepackage{amsthm}
\usepackage{amsmath}
\usepackage{graphicx}
\usepackage{amssymb}
\usepackage{color}

\newcommand{\la}{\lambda}
\newcommand{\vp}{{\mathbf p}}

\newcommand{\vz}{{\mathbf z}}

\newcommand{\vf}{{\mathbf f}}

\usepackage{tikz}
\usetikzlibrary{arrows}
\usepackage{float}

\begin{document}

\title{\textbf{Upper bound on the rate of convergence and truncation bound
for non-homogeneous birth and death processes on $\mathbb{Z}$}}

%Different institutions

\author{{Y. A. Satin\footnote{Vologda State University; e-mail yacovi@mail.ru}},
{R. V. Razumchik\footnote{Federal Research Center ``Computer Science and Control'' of the Russian Academy of Sciences; Moscow Center for Fundamental and Applied Mathematics; e-mail rrazumchik@ipiran.ru}},
{A. I. Zeifman\footnote{{Vologda State University;  Federal Research Center ``Computer Science and Control'' of the Russian Academy of Sciences; Moscow Center for Fundamental and Applied Mathematics; Vologda Research Center RAS; e-mail a$\_$zeifman@mail.ru}}},
{I. A. Kovalev\footnote{Vologda State University; e-mail kovalev.iv96@yandex.ru}}}

\date{} %leave it as it is, do NOT remove

\maketitle

%% Do not use the following keywords and abstract environment, use \section instead to get desired layout
%%\begin{keywords}
%%Agent-based modelling; Genetic algorithms
%%\end{keywords}
%%
%%\begin{abstract}
%%Abstract goes here. This is an example of using the ECMS \LaTeX\ class.
%%\end{abstract}

%% Notes about sections: Use star * to avoid numbering of sections

%\section*{\textbf{KEYWORDS}}

%Inhomogeneous birth-death processes; queueing models; two-sided uniform approximation bounds

\section*{\textbf{Abstract}}

%We study inhomogeneous random walk on Z.
%We consider a class of inhomogeneous birth-death queueing models and  obtain uniform approximation bounds of two-sided truncations. Some examples are studied.

We consider the well-known problem of the computation of the (limiting) time-dependent performance
characteristics of one-dimensional continuous-time birth and death processes 
on $\mathbb{Z}$ with time varying and possible state-dependent intensities.
First in the literature upper bounds on the rate of convergence
along with one new concentration inequality are provided.
Upper bounds for the error of truncation are also given.
Condition under which a limiting (time-dependent) distribution exists
is formulated but relies on the quantities that need to be guessed
in each use-case. The developed theory is illustrated by two numerical examples within the queueing theory context.

\section{\textbf{Introduction}}

In this paper consideration is given to the random walk on the integers,
performed by a~particle, which takes only unit steps 
either to the left or to the right.
Its initial position may be arbitrary but fixed.
The main quantity under the consideration is
the position $X(t) \in \mathbb{Z}$ of the particle at time~$t$.
Yet meaningful statements related to its average position
$\mathsf{E}X(t)$ given that initially it was
at the origin (to be understood here as ${X(0)=0}$) 
will also be given. The particle's position $X(t)$ at time~$t$ is governed by the
two independent Poisson processes with possible time-dependent and state-dependent
parameters; henceforth if ${X(t)=i}$ at some time~$t$ then $\lambda_i(t)$ and $\mu_i(t)$ denote
the motion intensities to the right and left respectively.
From the other point of view the $X(t)$ can be viewed as the 
non-homogeneous birth and death process (BDP) on~$\mathbb{Z}$
--- a model used for numerous problem instances in finance,
genetics, biology, chemistry, physics etc.
Just for an example one can refer to the bibliography (up to 1982)
in~\cite{introx1}, which contains more than 300 papers
on the use of~BDP in the latter two subjects;
a more recent review (up to 2004) can be found in~\cite{introx1new}.
One intuitively clear example of $X(t)$ (which will be revisited
further in the numerical section) is provided by one problem known
in the literature as the taxicab problem~\cite{introx2,introx3}.
There is one queueing point whereto both taxis and passengers arrive 
one by one in~accordance with the two independent Poisson flows
possibly with time-dependent and possibly state-dependent arrival intensities.
The queue length may take any integer value:
negative values mean that there are passengers waiting
for taxis, whereas positive values mean that there are taxis 
waiting for passengers. Whenever the queue length 
is zero, the queueing point is free from both passengers and taxis.
From the given description is can be seen 
that $X(t)$ can be represented as the difference between 
the two Poisson variables.
If~the intensities depend on the
state of $X(t)$, it implies that
the admission of passengers/taxis
to the queueing-point is dynamically controlled.
Since the seminal paper \cite{introx2}
such queues and similar to $X(t)$ processes
have been the subject of extensive research and now
they are usually referred to as
double-sided or double-ended queues
(see \cite{Wang2021,introx4}),
unrestricted random walks on lattice
and bilateral BDPs \cite{introx5,introx6,introx7}.
Another intuitive but otherwise artificial example
(which will also be revisited in the numerical section)
is the system size/queue length in common~queueing systems\footnote{Not arbitrary ones, but only 
those in which the queue-size may change by at most~1 at a~time.} at epoch~$t$.
If one removes the impenetrable barrier 
at the origin, which means that the departures are also allowed, 
when the system size is zero or negative,
one arrives at another instance of~$X(t)$ (see~\cite{ex2,introx9}).
At last, another example can be extracted from the Markov
predator-prey models (or other models of species coexistence \cite{Zeifman1982}),
in which $X(t)$ is the difference between the predator and prey populations.

Bilateral non-homogeneous BDPs like $X(t)$ have already
been analyzed in the literature from various perspectives;
see, for example, \cite[Section~1]{Giorno2019}.
The basic questions under consideration are:
the computation of the time-dependent and the limiting probability distribution,
methods for the approximation of their transient behaviour,
determination of first-passage time densities via analytical and numerical methods.
The literature review, which we have been able to make,
shows that for one of the general cases
i.e. when the state space of $X(t)$ is $\mathbb{Z}$
and its transition intensities are allowed 
to be time- and state-dependent, most of the questions remain open.
The only feasible way to deal with such~$X(t)$ 
seems to be extensive use of numerical schemes 
for systems of ordinary differential equations (ODEs).
For the numerical approaches to be efficient,
in the first place one needs to know how to determine a~priori the
points of convergence and, in the cases when the ODE system is 
infinite, how to choose the truncation thresholds.
In this paper we show that the technique utilizing the 
notion of the logarithmic norm
and already available for the BDPs on the non-negative integers,
can be generalized for the BDPs on~$\mathbb{Z}$.
The theoretical results which follow are applicable only
to those cases when the limiting ergodic distribution
exists. The sufficient condition for that is  being formulated
(see \textit{Theorem~1}).

%But is complex and has to be checked for each problem instance.
%We note that the intensities are not going to zero.

%The typical questions of interest is the
%time-dependent distribution

%Since then the body of literature devoted to various BDPs
%grew so big, that a thorough review of
%the current state-of-art in the field requires a~separate paper.

%, where some results for
%one special case of the considered problem are obtained.
%Due to the reasons which will become clear shortly,
%those results are not covered by the theory presented below.
%possible time-dependent and state-dependent parameters;

%(henceforth --- random walk on $\mathbb{Z}$)
%The the intensity parameter is allowed to  time-dependent
%and statemay depend on time.

The purpose of this paper is two-fold.
Firstly we derive first in the literature
explicit upper bounds for the rate of convergence
of non-homogeneous BDPs on $\mathbb{Z}$
to~the limiting regime (whenever it exists).
The class of processes considered includes those
with all the transition intensities being
possibly time-varying and state-dependent, but bounded (see~\eqref{ub0}).
Secondly, we derive truncation bounds
(see \textit{Theorem~2}), which allow one to obtain numerical solutions
with the desired accuracy.
By virtue of two numerical experiments
it is demonstrated that this result may be particularly
useful for obtaining the limiting values of 
the time-dependent probabilities.

The questions of convergence of non-homogeneous BDPs
(and especially homogeneous) have been considered in 
many research papers.
The approach used here to obtain the results related to the convergence and 
truncation bounds is, of course, not new. 
It is based on the theory developed in the series of papers
by the authors. Basically it relies on the well-known 
connection between the transition
matrix of a~Markov chain and the corresponding
ODEs (specifically, Kolmogorov's forward equations).
The main ingredient is the norion of the logarithmic norm
of an~operator function and those estimates 
for the differential equations, which are available in 
the literature.
Using this approach in the previous papers it
was possible to obtain explicit
upper bounds for the
distance between two probability distributions 
(in some special norms) of the BDPs with either finite or 
countable (in one direction) state space i.e. $\mathbb{Z}^+$.
Here we show, that the approach can be generalized
to deal with quite general BDPs on the whole set~$\mathbb{Z}$.
Surprisingly this generalization does not come at price:
the upper bounds obtained for the case of $\mathbb{Z}$
are not weaker that in the case of $\mathbb{Z}^+$. 

%not lower in case of~$\mathbb{Z}$)

In what follows $\|\cdot\|$ denotes the $l_1$-norm, i.e.
if ${\bf x}$ is a column vector then
$\|{{\bf x}}\|=\sum_{k} |x_k|$.
Clearly, $\|{{\bf x}}\|=1$ if ${\bf x}$ is a probability vector.
The operator norm is assumed to be induced by the $l_1$-norm
on column vectors i.e. for any~linear operator~$A$ we have
$\|A \| = \sup_{j} \sum_{i} |a_{ij}|$.

%Note that the linear map given
%by the block matrix $(D^*)^{-1}={\footnotesize \begin{pmatrix} (D^*_L)^{-1} & 0 \\ 0 & (D^*_U)^{-1} \end{pmatrix}}$,

\section{\textbf{Preliminaries}}

Let $\{X(t), t \ge 0\}$ be the BDP with the state
space $\mathbb{Z}$ and 
the generators $\{Q(t)=(q_{ij}(t)), \ t \ge 0 \}$ defined by
$$
q_{i,i+1}\left( t\right) =\lambda_i(t),
\
q_{i,i-1}\left( t\right) = \mu_i(t)
\mbox{ and }
q_{ii}(t) = - (\lambda_i(t)+\mu_i(t)),
$$

\noindent In what follows  $\lambda_i(t)$ and $\mu_i(t)$ are
assumed to be non-random locally integrable for $t \in [0,\infty)$
continuous functions, satisfying
\begin{equation}
\label{ub0}
0 \le \lambda_i(t) \le {\overline \lambda}_i \le \Delta < \infty, \ \
0 \le  \mu_i(t) \le
{\overline \mu}_i \le \Delta < \infty
\end{equation}

\noindent for all $t \ge 0$,  $i \in \mathbb{Z}$ and some
constants ${\{{\overline \lambda}_i, i \in \mathbb{Z} \}}$,
${\{{\overline \mu}_i, i\in \mathbb{Z} \}}$ and~${\Delta}$.
The transition diagram of $X(t)$ is shown in the figure below\footnote{In order to keep the figure and the matrices readable, whenever it does not introduce any ambiguity,
the argument~$t$ of the intensity functions is omitted.}.

\begin{figure}[h]
	\centerline{
		\includegraphics[width=300pt]{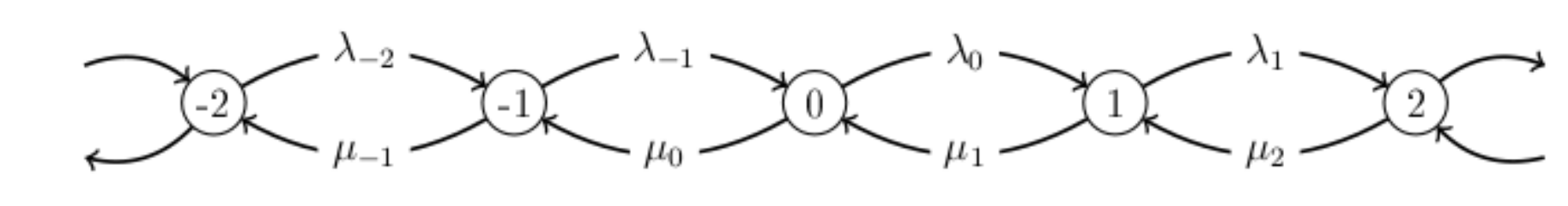}}
	\caption{Possible transitions for $X(t)$ and corresponding intensities}
\end{figure}

Let $p_{i}(t)=P\left\{ X(t)=i\right\}$ and
${\bf p}(t) = \left(\dots, {p}_{\text{-}1}(t), p_{0}(t), p_{1}(t), \dots \right)^T$.
For what follows it will be convenient to
write the Kolmogorov forward equations for the distribution of~$X(t)$ as\footnote{Special
cases of these well-known equations have been the starting point
for numerous papers; see, for example, early research on population dynamics in~\cite[Section~3]{introx8}.}
\begin{equation}
\label{ur01}
\frac{d}{dt}{\bf p}(t)=A(t){\bf p}(t),
\ t \ge 0,
\end{equation}

%in the infinite-dimensional vector space with the norm $l_1$

\noindent where $A(t)=(a_{ij}(t))$ is the transposed generator
i.e. $a_{ij}(t)=q_{ji}(t)$.
Since ${\|A(t)\| = 2\sup_{i \in \mathbb{Z}}(\lambda_{i}(t)+ \mu_i(t)) \le 4\Delta}$,
the linear operator $A(t)$ is bounded
and locally integrable for $t\in [0,\infty)$.
Thus (\ref{ur01}) is the system of differential
equations in the~space $l_1$ with the (bounded) linear operator.
And thus (see, for instance, \cite{DK})
it has the unique solution for
arbitrary initial conditions.
Moreover, if for some $s \ge 0$ the
probabilities ${p}_{i}(s)$ are all non-negative
for ${i \in \mathbb{Z}}$ and $\|{\bf p}(s)\|=1$,
then the same holds for ${\bf p}(t)$ when ${t \ge s}$.
What follows next relies on the concept of the logarithmic norm of
locally integrable operator functions
(see \cite{DK,newlo}) and available from the literature
estimates for differential equations; detailed definitions and derivations
can be recovered from, for example, \cite[Appendix]{z95}.

%Let the integer-valued time-dependent random variable $X(t)$ denote the
%total number of customers in the system at time $t \ge 0$.
%Then $X(t)$ is the CTMC with the state space $\{ 0, 1, 2 \dots \}$.
%Its transposed time-dependent  intensity matrix (generator)
%$A(t)=(a_{ij}(t))_{i,j=0}^{\infty}$ has the form

%\bigskip

%\section{\textbf{Upper bounds on the rate of convergence}}

\section{\textbf{Basic estimates}}

%В этом разделе будет доказано несколько теорем о
%нижних и верхних оценках для скорости сходимости к
%предельному режиму, а также
%получены новые оценки для вероятностей состояний,
%справедливые для любого $t \ge 0$.
%Начнем с теоремы о верхней оценке.

%${\bf Z}$--indexed series
%a not necessarily positive function $\beta^{**}(t)$
%such that for any ${k \in {\bf Z}\setminus \{0\}}$ and for almost all ${t\ge0}$}

\textit{Theorem 1. Let there exist a doubly infinite sequence of positive numbers $\{d_k, k= \pm 1, \pm 2 ,\dots\}$
such that ${\inf_{k \in \mathbb{Z} \setminus \{0\}} d_k=1}$ and ${\int_0^\infty \beta^{**}(u)du=\infty}$,
where ${\beta^{**}(t)=\inf_{k \in \mathbb{Z} \setminus \{0\}} \beta^{**}_k(t)}$
and the function ${\beta^{**}_k(t)}$ is given by }
$$
\beta^{**}_k(t)=
\begin{cases} \lambda _k\left( t\right) +\mu
_{k+1}\left(t\right) - \frac{d_{k+1}}{d_k} \lambda
_{k+1}\left(t\right) -  \frac{d_{k-1}}{d_k} \mu
_k\left( t\right), & k < -1\\
\lambda _{-1}\left( t\right) +\mu _{0}\left(t\right) -
\frac{d_{1}}{d_{-1}} \lambda _{0}\left(t\right) -
\frac{d_{-2}}{d_{-1}} \mu
_{-1}\left( t\right), & k = -1\\
\lambda _0\left( t\right) +\mu_{1}\left( t\right) -
\frac{d_{2}}{d_{1}} \lambda_{1}\left( t\right) -
\frac{d_{-1}}{d_{1}} \mu
_0\left( t\right), & k= 1,\\
\lambda _{k-1}\left( t\right) +\mu _{k}\left( t\right) -
\frac{d_{k+1}}{d_{k}} \lambda _{k}\left( t\right) -
\frac{d_{k-1}}{d_{k}} \mu
_{k-1}\left( t\right), & k > 1.
\end{cases}
.
$$
\textit{Then $X(t)$ is weakly ergodic and
for all $t\ge0$ and any initial conditions ${\bf p}^{*}(0)$ and $ {\bf p}^{**}(0)$
it holds that}
\begin{eqnarray}
\left\|{\bf p}^{*}(t)-{\bf p}^{**}(t)\right\|
\le
e^{-\int\limits_0^t \beta^{**}(u)\, du}
\sum\limits_{k \in \mathbb{Z} \setminus \{0\}}^N
(p^{*}_k(0)-p^{**}_k(0) )
\sum_{j=\min(1,k)}^{\max(-1,k)} d_j.
%\left\|{\bf p}^{*}(0)-{\bf p}^{**}(0)\right\|.
%\sum\limits_{\substack{k=-\infty \\ k \neq 0}}^\infty d_k |p^*_k(0)-p^{**}_k(0)|.
\label{the1}
\end{eqnarray}

%$p_0(t) = 1 - \sum\limits_{\substack{k=-\infty \\ k \neq 0}}^\infty p_k(t)$,

\textit{Proof.}
Since $p_0(t) = 1 - \sum_{k \in \mathbb{Z} \setminus \{0\}}^\infty p_k(t)$,
the Kolmogorov forward equations  (\ref{ur01}) for the distribution of $X(t)$
can be re-written as
\begin{equation}
\label{216}
\frac{d {}}{dt} \vz(t)= B(t){\vz}(t)+{\vf}(t),
\end{equation}
\noindent where the vectors ${\vf}(t)$ and ${\vz}(t)$ are
$$
{\bf f}\left(t\right)=\left(\dots,0, \mu_0(t), \la_0(t), 0, \dots
\right)^T, \ \
{\vz}(t) = \left(\dots, {p}_{\text{-}2}(t),{p}_{\text{-}1}(t), p_{1}(t), p_{2}(t), \dots \right)^T,
$$

%\footnote{Whenever
%it does not introduce any ambiguity, the argument $t$ of the intensity functions is %omitted for brevity.}

\noindent and the linear transformation $B(t)$ is given by the block matrix
$$
B(t)={\footnotesize \begin{pmatrix} B_{11}(t) & B_{12}(t) \\ B_{21}(t) & B_{22}(t) \end{pmatrix}},
$$
which entries $B_{ij}(t)$ are itself matrices of the form
{\footnotesize
$$
B_{11}(t)=\left(
\begin{array}{ccccccc}
\ddots & \vdots & \vdots & \vdots \\
 \cdots & - \left ({\la}_{\text{-}3}\!+\!{\mu}_{\text{-}3} \right )   & {\mu}_{\text{-}2} & 0  \\
 \cdots & {\la}_{\text{-}3}  & - \left ({\la}_{\text{-}2}\!+\!{\mu}_{\text{-}2} \right ) & {\mu}_{\text{-}1}   \\
 \cdots & - {\mu}_{0}  & {\la}_{\text{-}2}\!-\!{\mu}_{0}  & -\left ({\la}_{\text{-}1}\!+\!{\mu}_{\text{-}1} \!+\! {\mu}_{0} \right )  \\
\end{array}
\right),
\ \
B_{12}(t)=\left(
\begin{array}{ccccccc}
\vdots & \vdots & \vdots &  \iddots \\
0   & 0 & 0 & \cdots \\
0   & 0  & 0  &  \cdots  \\
 -{\mu}_{0}   &  -{\mu}_{0}  &  -{\mu}_{0} &   \cdots
\end{array}
\right),
$$
$$
B_{21}(t)=\left(
\begin{array}{ccccccc}
\cdots & -\la_{0}   & -\la_{0}  &  -\la_{0} \\
\cdots &0   & 0 & 0  \\
\cdots &0   & 0  & 0   \\
 \iddots & \vdots & \vdots & \vdots \\
\end{array}
\right),
\ \
B_{22}(t)=
%\begingroup % keep the change local
%\setlength\arraycolsep{1.5pt}
\left(
\begin{array}{ccccccc}
 - \left ({\lambda}_{1} \!+\! {\mu}_{1} \!+\! {\lambda}_{0} \right )    & {\mu}_{2} \!-\! {\lambda}_{0} & -\lambda_{0} & \cdots \\
{\lambda}_{1}   & -({\mu}_{2} \!+\! {\lambda}_{2})   & \mu_3  &   \cdots  \\
0  & {\lambda}_{2}  & -({\mu}_{3} \!+\! {\lambda}_{3})   &   \cdots \\
\vdots & \vdots & \vdots & \ddots \\
\end{array}
\right)
%\endgroup.
$$
%\begin{equation*}
%\label{1.01}
%\end{equation*}
%\begin{equation*}
%B(t)=
%\begingroup % keep the change local
%\setlength\arraycolsep{1.5pt}
%\begin{pmatrix}
%\ddots & \ddots & \ddots &   &  &   \cr
%%%%%%
%\ddots & - \left (\bar{\la}_{2}(t)+\bar{\mu}_{2}(t) \right )  & \bar{\mu}_{1}(t) &  0 & 0 &   \cr
%%%%%%
%\ddots & \left (\bar{\la}_{2}(t)-{\mu}_{0}(t) \right )  & - \left (\bar{\la}_{1}(t)+\bar{\mu}_{1}(t) +{\mu}_{0}(t) \right ) &  -{\mu}_{0}(t)  & - {\mu}_{0}(t) &   \cr
%%%%%%
% &  -\la_{0}(t)  & - {\lambda}_{0}(t) &  - \left ({\lambda}_{1}(t)+{\mu}_{1}(t) +{\lambda}_{0}(t) \right )  & \left ({\mu}_{2}(t) -{\lambda}_{0}(t) \right ) &  \ddots \cr
%%%%%%%%%%%%%%%%%%
% & 0  & 0 &  {\lambda}_{1}(t) & - \left ({\mu}_{2}(t)+{\lambda}_{2}(t) \right ) &  \ddots \cr
% &  &  &  \ddots & \ddots &  \ddots \cr
%\end{pmatrix} .
%\endgroup
%%\nonumber
%\label{1.01}
%\end{equation*}
}

\noindent Denote by $D^*_U$ and $D^*_L$
correspondingly the upper and
the lower triangular matrix of the form

{\footnotesize
$$
D^*_U=\left(
\begin{array}{ccccccc}
1   & 1 & 1 & \cdots \\
0   & 1  & 1  &   \cdots  \\
0   & 0  & 1  &   \cdots \\
\vdots & \vdots & \vdots & \ddots \\
\end{array}
\right),
\ \
D^*_L=\left(
\begin{array}{ccccccc}
\ddots & \vdots & \vdots & \vdots \\
 \cdots & - 1   & 0 & 0  \\
 \cdots & - 1   & -1 & 0  \\
 \cdots & - 1   & -1 & -1  \\
\end{array}
\right),
$$
}

\noindent Both of these matrices are known as semicirculant matrices.
Consider the linear transformation
given by the block matrix $D^*={\footnotesize \begin{pmatrix} D^*_L & 0 \\ 0 & D^*_U \end{pmatrix}}$.
In what follows we will need the
inverse linear map of $D^*$, which is further denoted by
$(D^*)^{-1}$. In order to show that it exists
for the considered matrix $D^*$ we will make use of the well-known
fact that the mapping of formal power series into the set
of infinite semicirculant matrices is an isomorphism.
Let us associate with the matrix $D^*_L$
the formal power series $P_L(z)=\sum_{i=0}^\infty z^i a_i$
(we write $P_L(z) \rightarrow D^*_L$).
The values of $a_i$ are in the first bottom row of $D^*_L$.
With the matrix $D^*_U$ we associate
the formal power series $P_U(z)=\sum_{i=0}^\infty z^i b_i$
(i.e. $P_U(z) \rightarrow D^*_U$).
The values of $b_i$ are in the first upper row of $D^*_U$.
Consider the matrix $P(z)={\footnotesize \begin{pmatrix} P_L(z) & 0 \\ 0 & P_U(z) \end{pmatrix}}$. Since the mapping $\rightarrow$ is an isomorphism,
then $P_L(z) \rightarrow D^*_L$ and $P_U(z) \rightarrow D^*_U$,
and thus then $P(z) \rightarrow D^*$.

%%%%%%%%%%%%%%%
%If the reciprocal series $(P_L(z))^{-1}$
%and $(P_U(z))^{-1}$ exist, then
%the inverse matrix is equal to

Note now that inverse matrix to $P(z)$, denote it by $(P(z))^{-1}$,
exists since $P_L(z) P_U(z) \neq 0$.
Denote the formal power series of ${1 \over P_L(z) P_U(z)}$
by $(P_L(z) P_U(z))^{-1}$.
Since $P_L(z) P_U(z)=\sum_{i=0}^\infty z^i (i+1)$,
then $(P_L(z) P_U(z))^{-1}=-1+2z-z^2$ (see \cite[Theorem~1.2b]{HEN}).
It is straightforward to check
that $(P_L(z) P_U(z))^{-1} P_U(z)=-1+z$
and $(P_L(z) P_U(z))^{-1} P_L(z)=1-z$.
Thus we have
\begin{multline}
(P(z))^{-1}
=
{1 \over P_L(z) P_U(z)}
{\footnotesize \begin{pmatrix} P_U(z) & 0 \\ 0 &  P_L(z)
 \end{pmatrix}}
 =
 \\
 =
 {\footnotesize \begin{pmatrix} (P_L(z) P_U(z))^{-1} P_U(z) & 0 \\ 0 &  (P_L(z) P_U(z))^{-1} P_L(z)
 \end{pmatrix}}
 =
 \\
 =
  {\footnotesize \begin{pmatrix} -1+z & 0 \\ 0 &  1-z
 \end{pmatrix}}.
\end{multline}

Both formal power series $-1+z$ and $1-z$ have associated semicitculant matrices
{\footnotesize
$$
1-z
\rightarrow
(D^*_U)^{-1} =\left(
\begin{array}{ccccccc}
1   & -1 & 0 & \cdots \\
0   & 1  & -1  &   \cdots  \\
0   & 0  & 1  &   \cdots \\
\vdots & \vdots & \vdots & \ddots \\
\end{array}
\right),
\ \
-1+z
\rightarrow
(D^*_L)^{-1}=\left(
\begin{array}{ccccccc}
\ddots & \vdots & \vdots & \vdots \\
 \cdots & - 1   & 0 & 0  \\
 \cdots &  1   & -1 & 0  \\
 \cdots & 0   & 1 & -1  \\
\end{array}
\right).
$$
}

Introduce the block matrix $(D^*)^{-1}={\footnotesize \begin{pmatrix} (D^*_L)^{-1} & 0 \\ 0 & (D^*_U)^{-1} \end{pmatrix}}$.
Thus we have $(P(z))^{-1}  \rightarrow (D^*)^{-1}$.
But since $P(z) (P(z))^{-1}=I$,
then $D^* (D^*)^{-1}=(D^*)^{-1}D^*=I$,
where $I$ is the identity matrix.
Thus $(D^*)^{-1}$ is the left and right inverse linear map of $D^*$.

Consider the similarity transformation $D^* B (t) (D^{*})^{-1}$,
further denoted by $B^*(t)$. It is well-defined and given by the matrix

%{\footnotesize
\begin{equation*}
B^*(t)=
\begingroup % keep the change local
\setlength\arraycolsep{1.5pt}
%\begin{pmatrix}
\left(
\begin{array}{ccc|ccc}
\ddots & \ddots &  &   &  &   \\
%%%%%
\ddots & - \left ({\la}_{\text{-}2}\!+\!{\mu}_{\text{-}1} \right )  & {\mu}_{\text{-}1} &  0 & 0 &   \\
%%%%%
 & {\la}_{\text{-}1}  & -\left ({\la}_{\text{-}1} \!+\! {\mu}_{0} \right ) &  {\mu}_{0}  & 0 &    \\ \hline
%%%%%
 & 0 & {\lambda}_{0} &  - \left ({\lambda}_{0} \!+\! {\mu}_{1} \right )  & {\mu}_{1} &   \\
%%%%%%%%%%%%%%%%%
 & 0  & 0 &  {\lambda}_{1} & - \left ({\mu}_{2} \!+\! {\lambda}_{1} \right ) &  \ddots \\
 &  &  &   & \ddots &  \ddots
%\end{pmatrix}
\end{array}
\right).
\endgroup
%\nonumber
\label{1.01}
\end{equation*}
%}

%{\color{blue} ?????надо бы проверить вид матрицы???????}

%$\prod\limits_{\substack{k=-\infty \\ k \neq 0}}^\infty d_k \neq 0$

\noindent Note that unlike the matrix $B(t)$
all off-diagonal entries of $B^*(t)$ are non-negative.
Choose an double infinite sequence $\{d_k, k= \pm 1, \pm 2 ,\dots\}$
of positive numbers and
consider the linear transformation ${D^{**}=diag\left(\dots,d_{-2},d_{-1},d_{1},d_{2},\dots\right)}$.
It is known (see \cite[p.~19]{cooke}) that $D^{**}$
has a unique right-hand reciprocal,
which is the diagonal matrix ${diag\left(\dots,1/d_{-2},1/d_{-1},1/d_{1},1/d_{2},\dots\right)0=(D^{**})^{-1}}$.
It is straightforward to check, that the
linear transformation $B^{**}(t)=D^{**} B^*(t) (D^{**})^{-1}$
is given by the matrix
%{\footnotesize
\begin{equation*}
B^{**}(t)=
\begingroup % keep the change local
\setlength\arraycolsep{1.5pt}
%\begin{pmatrix}
\left(
\begin{array}{ccc|ccc}
\ddots & \ddots &  &   &  &   \cr
%%%%%
\ddots & - \left ({\la}_{\text{-}2}\!+\!{\mu}_{\text{-}1} \right )  &  \frac{d_{\text{-}2}}{d_{\text{-}1}} {\mu}_{\text{-}1} &  0 & 0 &   \cr
%%%%%
 & \frac{d_{\text{-}1}}{d_{\text{-}2}}{\la}_{\text{-}1}  & -\left ({\la}_{\text{-}1} \!+\! {\mu}_{0} \right ) &   \frac{d_{\text{-}1}}{d_{1}}{\mu}_{0}  & 0 &   \\ \hline
%%%%%
 & 0 & {\lambda}_{0}\frac{d_{1}}{d_{\text{-}1}} &  - \left ({\lambda}_{0} \!+\! {\mu}_{1} \right )  & \frac{d_{1}}{d_{2}}{\mu}_{1} &   \\
%%%%%%%%%%%%%%%%%
 & 0  & 0 &  \frac{d_{2}}{d_{1}} {\lambda}_{1} & - \left ({\mu}_{2} \!+\! {\lambda}_{1} \right ) &  \ddots \\
 &  &  &   & \ddots &  \ddots
%\end{pmatrix} ,
\end{array}
\right),
\endgroup
%\nonumber
\label{1.01}
\end{equation*}
%}

% $D^{**} B^*(t) D^{**}^{-1}$
\noindent which has only non-negative off-diagonal elements.

Coming back to \eqref{216}, note that any
upper bound on the convergence rate  
to the limiting regime for $X(t)$,
corresponds to the same bound for 
the solutions of the system
\begin{equation}
\label{216n}
\frac{d {}}{dt} {\bf y} (t)= B(t){\bf y}(t),
\end{equation}
\noindent without the free term $ {\bf f}(t)$.
Here the vector ${\bf y}(t) =  \left(\dots, {y}_{\text{-}2}(t),{y}_{\text{-}1}(t), y_{1}(t), y_{2}(t), \dots \right)^T$ and its elements can either positive or negative. Denote $D=D^{**} D^*$ and ${\bf u}(t)=D {\bf u}(t)$. By left-multiplying both parts of~\eqref{216n} by~$D$, we get
\begin{equation}
\label{216nn}
\frac{d}{dt} {\bf u}(t) = B^{**}(t) {\bf u}(t),
\end{equation}

\noindent where ${\bf u}(t) = \left(\dots, {u}_{\text{-}2}(t),{u}_{\text{-}1}(t), u_{1}(t), u_{2}(t), \dots \right)^T$ is, as well as ${\bf y}(t)$, the vector
with the elements of arbitrary signs. Let us estimate the logarithmic norm of $B^{**}(t)$.
%, which we further denote by $\gamma\left(\cdot\right)$.
It is well-known that in the $l_1$-norm
the logarithmic norm
of a  (locally integrable) operator $F(t)=(f_{ij}(t))$ is equal to
${\sup_i \left ( f_{ii}(t) + \sum_{j \neq i} |f_{ji}(t)| \right )=\gamma\left( F(t) \right)}$ (see, for example, \cite[Appendix]{z95}).
By direct inspection it can be instantly seen that
the $k$th column sum of $B^{**}(t)$ is equal to
$-\beta^{**}_k(t)$, where
\begin{eqnarray*}
\beta^{**}_{k}\left( t\right) =
\begin{cases} \lambda _k\left( t\right) +\mu
_{k+1}\left(t\right) - \frac{d_{k+1}}{d_k} \lambda
_{k+1}\left(t\right) -  \frac{d_{k-1}}{d_k} \mu
_k\left( t\right), & k < -1\\
\lambda _{-1}\left( t\right) +\mu _{0}\left(t\right) -
\frac{d_{1}}{d_{-1}} \lambda _{0}\left(t\right) -
\frac{d_{-2}}{d_{-1}} \mu
_{-1}\left( t\right), & k = -1\\
\lambda _0\left( t\right) +\mu_{1}\left( t\right) -
\frac{d_{2}}{d_{1}} \lambda_{1}\left( t\right) -
\frac{d_{-1}}{d_{1}} \mu
_0\left( t\right), & k= 1,\\
\lambda _{k-1}\left( t\right) +\mu _{k}\left( t\right) -
\frac{d_{k+1}}{d_{k}} \lambda _{k}\left( t\right) -
\frac{d_{k-1}}{d_{k}} \mu
_{k-1}\left( t\right), & k > 1.
\end{cases}
%}
 \label{211}
\end{eqnarray*}

\noindent Thus $\gamma\left( B^{**}(t) \right)=\sup\limits_{k \in \mathbb{Z} \setminus \{0\}} \, (- \beta^{**}_{k} (t))=
-\inf\limits_{k \in \mathbb{Z} \setminus \{0\}} \, \beta^{**}_{k} (t)$.
Let ${\{d_k, k= \pm 1, \pm 2 ,\dots\}}$ be such a doubly infinite sequence,
that ${\inf\limits_{k \in \mathbb{Z} \setminus \{0\}} \, \beta^{**}_{k} (t)<\infty}$
for ${t \ge 0}$. Denote ${\inf\limits_{k \in \mathbb{Z} \setminus \{0\}} \, \beta^{**}_{k} (t)=\beta^{**}(t)}$.
Then
$$
\| B^{**}(t) \| \le 4 \Delta - \beta^{**}(t)
$$

\noindent and thus $B^{**}(t)$ is the bounded operator.
Now, if $V(t,z)$ is the Cauchy operator of the equation (\ref{216nn}), then
for any~$t$ and~$s$
the following bound holds (for the justification see, for example,~\cite[Theorem~A2]{z95}):
\begin{equation}
\|V(t,s)\| \le e^{-\int\limits_s^t \beta^{**}(u) \, d u},
\ {0 \le s \le t}.
\label{220}
\end{equation}

\noindent Now let ${\bf p}^{*}(t)$ and ${\bf p}^{**}(t)$ be such that
the corresponding  $D{\bf z}^{*}(t)$ and $D{\bf z}^{**}(t)$ exist.
Then for any $t \ge 0$ we have
\begin{eqnarray*}
\left\|{\bf p}^{*}(t)-{\bf p}^{**}(t)\right\|
&& \le
2 \left\|{\bf z}^{*}(t)-{\bf z}^{**}(t)\right\|
\le
\\
&& \le
 \left\|D{\bf z}^{*}(t)-D{\bf z}^{**}(t)\right\|
\le
\\
&&\le
e^{-\int\limits_0^t \beta^{**}(u)\, du}\left\|
D{\bf z}^{*}(0)- D{\bf z}^{**}(0)\right\|.
\end{eqnarray*}

\begin{flushright}
$\qedsymbol$
\end{flushright}

The inequality \eqref{the1} holds even if $\int_0^\infty \beta^{**}(u)du<\infty$.
This happens only if the intensities approach~$0$ as~$t$ becomes infinite
and thus the limiting ergodic distribution cannot not~exist (cf. \cite[Example~3]{introx3}).
It is also worth noticing here that $\beta^{**}(t)$ is not necessarily an
everywhere positive function.

\bigskip

\textit{Corollary 1. Assume that under the assumptions of the Theorem 1
there exist positive constants $M$ and $\beta^{**}$ such that
$e^{-\int_s^t \beta^{**}(\tau)\, d\tau} \le M e^{-\beta^{**}\cdot (t-s)}$
for any ${0 \le s \le t}$.
Then for any positive integer $N$ and all $t \ge 0$ it holds that}
\begin{equation}
\label{corol1}
Pr(|X(t| \ge N) \le
M  \left (
  e^{-\beta^{**} t}
\sum\limits_{\substack{k=-N \\ k \neq 0}}^N  p_k(0)
\sum_{j=\min(1,k)}^{\max(-1,k)} d_j
+
\frac{d_{-1}\overline{\mu}_0 + d_{1} \overline{\lambda}_0}{\beta^{**}}
\right )
\left ( {
\sum\limits_{\substack{j=-N \\ j \neq 0}}^{N} d_j
\over \sum\limits_{j=-N}^{-1} d_j \cdot \sum\limits_{j=1}^{N} d_j } \right ).
\end{equation}

\textit{Proof.} Consider~\eqref{216} and note that its solution is 
\begin{equation}
\label{soldif}
{\bf z}(t)
= V(t,0){\bf z}(0)+\int_0^t{}V(t,\tau){\bf
f}(\tau)\,d\tau.
\end{equation}

\noindent Let us left-multiply the left and the right part
of the previous relation by~$D$. Using the estimates obtained 
in \textit{the Theorem~1} and assuming that
there exist constants $M>0$ and $\beta^{**}>0$ such that
 $e^{-\int\limits_s^t \beta^{**}(u)\, du} \le M e^{-\beta^{**} \cdot (t-s)}$
 for any ${0 \le s \le t}$, we get
\begin{multline}
\label{leftD1}
\| D {\bf z}\left(t\right)\| \le \|V\left(t,0\right)\|
\|D {\bf z}\left(0\right)\|
+
\int_0^t \|V\left(t,s\right)\|
\|D {\bf f}\left(s\right)\| \, ds \le
 \\
 \le
 M e^{-\beta^{**} t}
\|D {\bf z}\left(0\right)\|
+
 \int_0^t M e^{-\beta^{**}\cdot (t-s)}
 \|D{\bf f}\left(s\right)\| \, ds  \le
 \\
 \le
  M
  e^{-\beta^{**} t}
\left (
\sum\limits_{k=-N}^{-1}
p_k\left(t\right)
\sum\limits_{j=k}^{-1} d_j
+
\sum\limits_{k=1}^N
p_k\left(t\right)
\sum\limits_{j=1}^{k} d_j
  \right )
 +
M \frac{d_{-1}\overline{\mu}_0 + d_{1} \overline{\lambda}_0}{\beta^{**}}.
\end{multline}

\noindent Indeed, $\|D {\bf f}(t)\| =d_{-1}\mu_0(t) + d_{1}\lambda_0(t)\le \left(d_{-1} \overline{\mu}_0 + d_{1} \overline{\lambda}_0\right)$.
Note that since all $d_k$ are positive, then for
any positive integer $N$ we have
\begin{eqnarray*}
\| D {\bf z}\left(t\right)\|
&=&
\cdots+\left(d_{-2}+d_{-1}\right)p_{-2}(t)+
d_{-1}p_{-1}(t)+d_{1}p_{1}(t)+\left(d_{1}+d_{2}\right)p_{2}(t)+\cdots
=
\\
&=&
\sum\limits_{k=-\infty}^{-1}
p_k\left(t\right)
\sum\limits_{j=k}^{-1} d_j
+
\sum\limits_{k=1}^\infty p_k\left(t\right)
\sum\limits_{j=1}^{k} d_j
\ge
\\
&\ge &
\sum\limits_{k=-\infty}^{-N}
p_k\left(t\right)
\sum\limits_{j=k}^{-1} d_j
+
\sum\limits_{k=N}^\infty p_k\left(t\right)
\sum\limits_{j=1}^{k} d_j.
\end{eqnarray*}

\noindent
From here we get two concentration inequalities for $X(t)$,
which are valid for any integer $N>0$:
$$
\sum\limits_{k=-\infty}^{-N} p_k\left(t\right) \le \frac{\|D{\bf z}\left(t\right)\|}{\sum\limits_{j=-N}^{-1} d_j},
\ \
\sum\limits_{k=N}^\infty p_k\left(t\right) \le \frac{\|D{\bf z}\left(t\right)\|}{\sum\limits_{j=1}^{N} d_j}.
$$
Combining this with the upper bound for $\|D{\bf z}\left(t\right)\|$
we get for any positive integer $N$:
$$
Pr(X(t)\le-N) \le
M
\left (
  e^{-\beta^{**} t}
\sum\limits_{\substack{k=-N \\ k \neq 0}}^N p_k(0)
\sum_{j=\min(1,k)}^{\max(-1,k)} d_j
+
\frac{d_{-1}\overline{\mu}_0 + d_{1} \overline{\lambda}_0}{\beta^{**}}
\right )
\left (\sum\limits_{j=-N}^{-1} d_j\right )^{-1},
$$
$$
Pr(X(t) \ge N) \le
M \left (
   e^{-\beta^{**} t}
\sum\limits_{\substack{k=-N \\ k \neq 0}}^N p_k(0)
\sum_{j=\min(1,k)}^{\max(-1,k)} d_j
+
 \frac{ d_{-1}\overline{\mu}_0 + d_{1} \overline{\lambda}_0}{\beta^{**}}
 \right )
\left (\sum\limits_{j=1}^{N} d_j \right )^{-1}.
$$

\begin{flushright}
$\qedsymbol$
\end{flushright}

Note that since the bound in the \textit{Corollary~1}
is valid for any $t$, it is valid for the limiting probabilities
as well (if they exist).
One can simplify the bound \eqref{corol1} by fixing the initial condition.
For example, if $p_0(0)=1$, which implies that ${\bf z}\left(0\right)={\bf 0}$,
then the first term in the brackets in the right hand-side of
\eqref{corol1} is zero.

%These inequalities can be represented as a single one for the random
%variable $|X(t)|$.

%${\bf Z}$--indexed series
%{\color {blue} ??? что можно сказать о нижней оценке скорости сходимости? ??? }

%\section{\textbf{TRUNCATIONS}}

%{\color {blue} ??? Зачем нужна теорема 2? ???? }

% as it was seen initial condition $i=0$ is important. Let us assume that
%initially we are in i=1

\bigskip

\textit{Theorem 2. Let $X(t)$ be a BDP on~$\mathbb{Z}$
for which the Corollary~1 holds. Let $X^*(t)$ be its truncated version
with the state space $\{N_1,\dots,0,\dots,N_2\}$,
${N_1<0}$, ${N_2>0}$. If there exist
constants $\beta^*$, $M^*$ and $\{d^*_k, {k \in \mathbb{Z}\setminus \{0\}}\}$,
such that the Corollary~1 holds for $X^*(t)$, then
the following upper bound for the difference between the probability
distributions of $X(t)$ and $X^*(t)$
}
\begin{equation}
\|{\bf p}\left(t\right)-{\bf p}^*\left(t\right)\|\le
\frac{4
M\,M^*\left({\overline \mu}_{0}\,d^*_{-1}+{\overline \lambda}_{0}\,d^*_{1}\right)}{d \beta^*\,\beta^{**}}
\left( \frac{ \sum\limits_{j=N_1-1}^{-1} d_j {\overline \mu}_{N_1} }{ \sum\limits_{j=N_1}^{-1} d^*_j}+
\frac{\sum\limits_{j=1}^{N_2+1} d_j {\overline \lambda}_{N_2}}{  \sum\limits_{j=1}^{N_2} d^*_j } \right)
\label{teor2}
\end{equation}

\noindent \textit{holds for any ${t\ge 0}$ if $X(0)=X^*(0)=0$.}

%$N_1,N_1+1, \cdots, 0, \cdots, N_2$

\textit{Proof.} Consider the BDP $X^*(t)$ with the state space $\mathbb{Z}$
and the intensities $\la^*_{k}(t)=\la_{k}(t)$ if $N_1 \le k <N_2$,
and $\mu^*_{k}(t)=\mu_{k}(t)$ if $N_1 < k \le N_2$
and other intensities equal to zero.
Thus the linear operator $A^*\left(t\right)$
is still given by the bi-infinite matrix.
The Kolmogorov forward equations for the distribution of
the $X^*(t)$, being the truncated $X(t)$,
are
\begin{equation}\label{wOutSubstTwo}
\frac{d}{dt}{\bf p}^*(t)=A^*\left(t\right){\bf p}^*\left(t\right),
\end{equation}

\noindent where
$$
{\bf p}^*(t) =  \left(\dots,0,{p}_{N_1}(t),\dots, {p}_{\text{-}2}(t),{p}_{\text{-}1}(t), p_{1}(t), p_{2}(t), \dots, {p}_{N_2}(t), 0,\dots \right)^T.
$$
\noindent Since $p^*_0(t) = 1 - \sum_{j \neq 0} p^*_j(t)$, then
we obtain from (\ref{wOutSubstTwo})
\begin{equation}\label{withSubstTwo2}
\frac{d}{dt}{\bf z}^*(t)=B^*\left(t\right){\bf z}^*\left(t\right)+{\bf f}^*\left(t\right).
\end{equation}

\noindent Rewrite (\ref{withSubstTwo2}) in the  form:
\begin{equation}\label{withSubstTwo}
\frac{d}{dt}{\bf z}^*(t)=B\left(t\right){\bf z}^*\left(t\right)+\left(B^*\left(t\right)-B\left(t\right)\right){\bf z}^*\left(t\right)+{\bf f}^*\left(t\right).
\end{equation}

\noindent Then we have the following relations between the solutions of (\ref{216})
and (\ref{withSubstTwo}):
{
\begin{multline}
%\begin{split}
{\bf z}\left(t\right)-{\bf z}^*\left(t\right) = V\left(t,0\right)\left({\bf z}\left(0\right)-{\bf z}^*\left(0\right)\right)+\\+
 \int_0^t
V\left(t,s\right)\left(B\left(s\right)-B^*\left(s\right)\right){\bf
z}^*\left(s\right)\,ds+ \\+ \int_0^t V\left(t,s\right)\left({\bf
f}\left(s\right)-{\bf f}^*\left(s\right)\right)\, ds,
%\end{split}
\label{zmz}
\end{multline}
}

\noindent where
$$
\left(B\left(s\right)-B^*\left(s\right)\right){\bf
z}^*\left(s\right) =  \left(
 \cdots,  0,
\mu_{N_1} p^*_{N_1}, -\mu_{N_1} p^*_{N_1},  0,
\cdots,
 0,  -\la_{N_2}p^*_{N_2},  \la_{N_2}p^*_{N_2},  0, \cdots
\right)^T.
$$

%, where $N_1<0 < N_2.$

For simplicity we assume further, that
${\bf z}\left(0\right)={\bf z}^*\left(0\right)=0$
(i.e. $X(0)=X^*(0)=0$ with the probability~$1$)
Then ${\bf f}\left(s\right)={\bf f}^*\left(s\right)$ for any~$s$.
Next, it is clear that the first and the third terms in the~\eqref{zmz}
are equal to zero and the difference
between ${\bf z}\left(t\right)$ and ${\bf z}^*\left(t\right)$ is just
\begin{equation}\label{subZ}
{\bf z}\left(t\right)-{\bf z}^*\left(t\right) = \int_0^t V\left(t,s\right)\left(B\left(s\right)-B^*\left(s\right)\right){\bf z}^*\left(s\right) \, ds.
\end{equation}

Let  $\{d^*_k, k= \pm 1, \pm 2 ,\dots\}$
be an double infinite sequence
of positive numbers such that
there exist positive $M^*$ and $\alpha^*$ such that
\begin{equation}
e^{-\int\limits_s^t \beta^*(\tau)\, d\tau} \le M^* e^{-(t-s)\cdot \beta^*},
\label{22002}
\end{equation}
\noindent for any $0 \le s \le t$, where $\beta^*\left( t\right)= \inf_{k \in \mathbb{Z} \setminus \{0\}} \beta^*_k\left( t\right)$
and the functions $\beta^*_k\left( t\right)$ is given by
$$
\beta^*_k\left( t\right) =
\begin{cases} \lambda^*_k\left( t\right) +\mu^*
	_{k+1}\left(t\right) - \frac{d^*_{k+1}}{d^*_k} \lambda
	^*_{k+1}\left(t\right) -  \frac{d^*_{k-1}}{d^*_k} \mu^*
	_k\left( t\right), & k < -1,\\
	\lambda^*_{-1}\left( t\right) +\mu^*_{0}\left(t\right) -
	\frac{d^*_{1}}{d^*_{-1}} \lambda^* _{0}\left(t\right) -
	\frac{d^*_{-2}}{d^*_{-1}} \mu^*
	_{-1}\left( t\right), & k = -1,\\
	\lambda^* _0\left( t\right) +\mu^*_{1}\left( t\right) -
	\frac{d^*_{2}}{d^*_{1}} \lambda^*_{1}\left( t\right) -
	\frac{d^*_{-1}}{d^*_{1}} \mu^*
	_0\left( t\right), & k= 1,\\
	\lambda^*_{k-1}\left( t\right) +\mu^*_{k}\left( t\right) -
	\frac{d^*_{k+1}}{d^*_{k}} \lambda^*_{k}\left( t\right) -
	\frac{d^*_{k-1}}{d^*_{k}} \mu^*	_{k-1}\left( t\right), & k > 1.\\
	\end{cases}
$$

%Consider the linear transformation ${D^{*}=diag\left(\dots,d^*_{-2},d^*_{-1},d^*_{1},d^*_{2},\dots\right)}$.
%Then
By left-multiplying both parts of
\eqref{subZ} by the matrix $D$,
introduced abive,
 and using the estimates obtained above we get:
\begin{multline}\label{subBp}
%\begin{split}
\|
D\left(B\left(s\right)-B^*\left(s\right)\right){\bf
z}^*\left(s\right)\|
\le\left|
\sum\limits_{j=N_1-1}^{-1} d_j
+
\sum\limits_{j=N_1}^{-1} d_j
\right|\mu_{N_1}(s) p^*_{N_1}(s)+ \\ +\left|
\sum\limits_{j=1}^{N_2+1} d_j
+
\sum\limits_{j=1}^{N_1} d_j
\right|\la_{N_2}(s)p^*_{N_2}(s) \le
2\sum\limits_{j=N_1-1}^{-1} d_j {\overline \mu}_{N_1} p^*_{N_1}(s) + 2
\sum\limits_{j=1}^{N_2+1} d_j
{\overline \lambda}_{N_2}p^*_{N_2}(s).
%\end{split}
\end{multline}

%\begin{equation}\label{estimate_p_kTwo}
%p^*_k\left(t\right) \le \left \{ \begin{array}{cc}
%\frac{M^* \left({\overline \mu}_0 d^*_{-1}+ {\overline \lambda}_0 d^*_{1}\right)}{\beta^* \, \sum\limits_{j=k}^{-1} %d^*_j},  \, k < 0\\
%\frac{M^* \left({\overline \mu}_0 d^*_{-1}+ {\overline \lambda}_0 d^*_{1}\right)}{\beta^* \, \sum\limits_{j=1}^{k} %d^*_j}, \, k > 0\\
%\end{array}
%\right. .
%\end{equation}

Since ${p^*_k\left(t\right) \le \frac{M^*}{\beta^*}
\left({\overline \mu}_0 d^*_{-1}+ {\overline \lambda}_0 d^*_{1}\right)( \sum_{j=\min(1,k)}^{\max(-1,k)} d^*_j)^{-1}}$ for any $k\neq 0$, then
under the assumption that both BDPs start in the $0^{th}$ state, relations (\ref{subZ})
and (\ref{subBp})
imply the bound
$$
\|D({\bf z}\left(t\right)-{\bf z}^*\left(t\right))\|\le \frac{2
M\,M^*\left({\overline \mu}_{0}\,d^*_{-1}+{\overline \lambda}_{0}\,d^*_{1}\right)}{\beta^*\,\beta^{**}}
\cdot\left( \frac{ \sum\limits_{j=N_1-1}^{-1} d_j {\overline \mu}_{N_1} }{ \sum\limits_{j=N_1}^{-1} d^*_j}+
\frac{\sum\limits_{j=1}^{N_2+1} d_j {\overline \lambda}_{N_2}}{  \sum\limits_{j=1}^{N_2} d^*_j } \right).
$$

% and (\ref{estimate_p_kTwo})

Put $d=\min\left(d_{-1},d_{1}\right)$. The following sequence
of inequalities completes the proof:
\begin{multline*}
\|{\bf p}\left(t\right)-{\bf p}^*\left(t\right))\|\le
\cdots+|p_{-1}(t)-p^*_{-1}(t)|+
|p_{0}(t)-p^*_{0}(t)|+
 |p_{1}(t)-p^*_{1}(t)|+\cdots\le
\\
%%%%%%%%%%%%%%
\le
\cdots+\frac{d_{-1}+d_{-2}}{d_{-1}+d_{-2}}\, |p_{-2}(t)-p^*_{-2}(t)|
+
\frac{ d_{-1}}{d_{-1}}\, |p_{-1}(t)-p^*_{-1}(t)|+
\\
+
\frac{d_{1}}{d_{1}}\,
|p_{1}(t)-p^*_{1}(t)|+\frac{d_{1}+d_2}{d_{1}+d_2}\,
|p_{2}(t)-p^*_{2}(t)|\cdots \le
\\
%%%%%%%%%%%%%%%%%%%
\cdots+\frac{d_{-1}+d_{-2}}{d}\, |p_{-2}(t)-p^*_{-2}(t)|
+
\frac{ d_{-1}}{d}\, |p_{-1}(t)-p^*_{-1}(t)|+
\\
+
\frac{d_{1}}{d}\,
|p_{1}(t)-p^*_{1}(t)|+\frac{d_{1}+d_2}{d}\,
|p_{2}(t)-p^*_{2}(t)|\cdots \le {1 \over d}
\|D({\bf z}\left(t\right)-{\bf z}^*\left(t\right))\|.
\end{multline*}

%This completes the proof.

\begin{flushright}
$\qedsymbol$
\end{flushright}

The argumentation of the \textit{Theorem~2} allows one also to
obtain the upper bound for the truncation error, when computing
the average value $\mathsf{E}X(t)$ given that initially the process
was in the $0^{th}$ state. Let $W=\inf_{k \ge 1}
\left (
\frac{\sum\limits_{j=-k}^{-1} d_j}{k},\,  \frac{ \sum\limits_{j=1}^{k} d_j}{k}
\right )$.
Then
\begin{multline*}
 \sum_{k \in \mathbb{Z} \setminus \{0\}}^\infty k |p_{k}(t)-p^*_{k}(t)|
 =
\cdots+1 \cdot |p_{-1}(t)-p^*_{-1}(t)|+
0 \cdot |p_{0}(t)-p^*_{0}(t)|+
1 \cdot |p_{1}(t)-p^*_{1}(t)|+\cdots\le
\\
%%%%%%%%%%%%%%
\le
\cdots+2 \cdot \frac{d_{-1}+d_{-2}}{d_{-1}+d_{-2}}\, |p_{-2}(t)-p^*_{-2}(t)|
+
1 \cdot \frac{ d_{-1}}{d_{-1}}\, |p_{-1}(t)-p^*_{-1}(t)|+
\\
+
1\cdot \frac{d_{1}}{d_{1}}\,
|p_{1}(t)-p^*_{1}(t)|+
2 \cdot \frac{d_{1}+d_2}{d_{1}+d_2}\,
|p_{2}(t)-p^*_{2}(t)|\cdots \le
\\
%%%%%%%%%%%%%%%%%%%
\cdots+\frac{d_{-1}+d_{-2}}{W}\, |p_{-2}(t)-p^*_{-2}(t)|
+
\frac{ d_{-1}}{W}\, |p_{-1}(t)-p^*_{-1}(t)|+
\\
+
\frac{d_{1}}{W}\,
|p_{1}(t)-p^*_{1}(t)|+\frac{d_{1}+d_2}{W}\,
|p_{2}(t)-p^*_{2}(t)|\cdots \le
 {1 \over W}
\|D({\bf z}\left(t\right)-{\bf z}^*\left(t\right))\|.
\end{multline*}

\noindent Using now the upper bound for
$\|D({\bf z}\left(t\right)-{\bf z}^*\left(t\right))\|$
from the \textit{Theorem~2},
we get the upper bound for the ${\sum_{k \in \mathbb{Z} \setminus \{0\}}^\infty k |p_{k}(t)-p^*_{k}(t)|}$.
In what follows we show,
how the developed theory can be used to
obtain explicit results.

\section{\textbf{Numerical examples}}

% {\color {blue} На мой взгляд, в примерах интересно было бы
% построить графики следующих функций: $E(X(t)|X(0)=0)$
% и $Var(X(t)|X(0)=0)$. Также интересны и индивидуальные
% вероятности, двумя различными способами.
% Первый способ:  по оси X отложены значения $t$,
% а по оси Y значения $P(X(t)=i)$. Можно выбрать несколько
% значений $i$ (скажем -10, -5, 0, 5, 10).
% Второй способ:
% по оси X отложены значения $i$ (скажем от -10 до 10),
% а по оси Y значения $P(X(t)=i)$. В качестве значений $t$ можно взять
% 1/8, 1/4,1/2,1,2,4,8,100. Или какие-то другие.}
%
%\subsection*{\textbf{Example 1}}

%the time intervals
%between successive steps to the right and to the left
%are independent and exponentially distributed with
%time-varying parameters.

Two examples are considered in this section.
Their main purpose is to illustrate that the developed theory
indeed allows one to study numerically
arbitrary bilateral BDP $X(t)$ with uniformly bounded
and state-dependent intensity functions.
Specific forms of the intensity functions
have been chosen for convenience of computation.
In each case it is assumed that ${X(0)=0}$.

% $X(t)$ starts at time ${t=0}$
%in the position ${X(0)=0}$.

In the first example we consider the randomized random walk on the integers,
say $X(t)$, which represents the position at time~$t$ of a
particle moving along, say $x$-axis, according to the following rules.
Its position $X(t)$ can be shifted by at most $1$ to the right or left,
and it is assumed that these changes are
governed by the two Poisson processes with the time-varying parameters. 
Specifically, when the particle is in a position
$i$ on the positive part of the $x$-axis, it will move to
the position $j$ in
the infinitesimal time $h>0$
with the probability
$$
Pr\left\{ X(t+h)=j\left| X(t)=i\right. \right\}
=
\begin{cases}
\lambda(t) h, & if j-i=1,\\
\mu_i(t) h, & j-i=1, \\
-(\lambda(t)+\mu_i(t)) h, & j=i,\\
0, & \mbox{ otherwise}.
\end{cases}
$$

\noindent The next position $j$
in  the infinitesimal time $h>0$
of the particle residing in the position $i$
on the negative part of the
$x$-axis is governed by the probability
$Pr\left\{ X(t+h)=-j\left| X(t)=-i\right. \right\}$.
If the particle enters the state $0$ then
its next state is $1$ or $-1$ with the probability $\lambda(t) h$.

We make further simplifications. Let us assume that
$\mu_i(t)=min (i; S) \mu(t)$.
Then when the particle is in the non-negative part
of the $x$-axis, then $X(t)$
represents the number of customers
present in the classic $M/M/S/\infty$ queue
at epoch $t$.
This example is somewhat artificial one and is due to~\cite{ex2}.

From the \textit{Theorem~1} one can obtain the upper bound
for the rate of convergence, if a~double infinite sequence,
say $\{d_k, k= \pm 1, \pm 2 ,\dots\}$, can be found such that
$\int_0^\infty \beta^{**}(u)du=\infty$.
Let us put ${d_1=1}$ and ${d_{k}=d^{k-1}}$ for ${k \ge 2}$, where ${d> 1}$.
Then we have:
\begin{eqnarray}
 \beta^{**}_{k}\left( t\right) =
\begin{cases}
\mu\left( t\right) -d \lambda\left( t\right), \ k= 1\\
\mu\left( t\right) - (d-1)\lambda\left( t\right) +\frac{(k-1)(d-1)}{d} \mu\left( t\right), \ 2 \le k \le S\\
\left(1-\frac{1}{d}\right)\left(S\mu\left( t\right) -d\lambda\left( t\right)\right), k > S.
\end{cases}
\label{2112}
\end{eqnarray}

\noindent Assume for now that there exists $\theta(t)$
such that $ \left(S\mu\left( t\right) -d\lambda\left( t\right)\right) \ge \theta(t) $.
Then $\beta^{**}\left( t\right) = \min \left(\mu\left( t\right) -d\lambda\left( t\right), \left(1-\frac{1}{d}\right)\theta(t) \right)$ and the upper bound follows from~\eqref{the1}.
Further insight can be gained if one fixes exact values of $S$, $\lambda_k(t)$ and $\mu_k(t)$.
So let us assume that ${S = 2}$, ${ \lambda\left( t\right) = 1+\sin(2\pi t)} $, $ {\mu_k\left( t\right)=3\min(k,S)}$. Then if one puts ${d = \frac{8}{7}}$ and ${d_{k}=\left(\frac{8}{7}\right)^{k-1}}$ for ${k \ge 1}$, then the constants $\beta^{**}$
and $M$ from the \textit{Theorem~1} and \textit{Corollary~1} are equal to
$\beta^{**}=\frac{13}{28}$ and $M=1$.
For the truncated process $X^*(t)$ with the truncation threshold $N=150$
one can put ${ d^* = \frac{4}{3}}$ and ${d^*_{k}=\left(\frac{4}{3}\right)^{k-1}}$ for ${k \ge 1}$. Then the constants $\beta^{*}$
and $M^*$ from the~\textit{Theorem~2} are equal to $\beta^*=\frac{1}{3}$ and $M^*=1$.
Thus, since  ${\|D(\vz(t)-\vz^*(t))\|\le 2 \times 10^{-8}}$
from~\eqref{teor2} one gets
$$
\|\vp(t)-\vp^*(t)\| \le 2 \times 10^{-8}.
$$
\noindent and from the comments, following the~\textit{Theorem~2},
one obtains
$$
\sum_{k \in \mathbb{Z} \setminus \{0\}}^\infty |k| |p_{k}(t)-p^*_{k}(t)|
\le 2 \times 10^{-8}.
$$

\begin{figure}[!ht]
	\centering
	\includegraphics[scale=0.4]{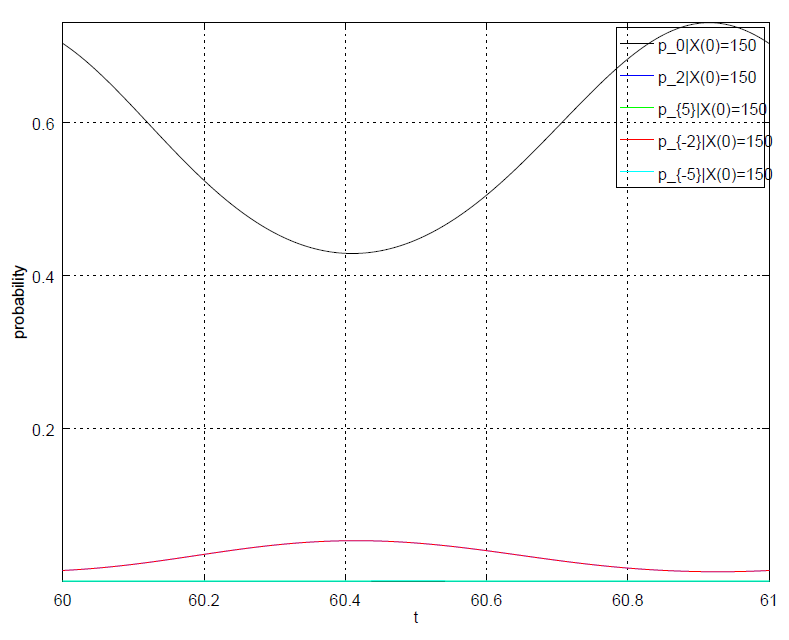}
	%\vspace{-1cm}
	\caption{Limiting probability of particle position $X(t)$ at time $t$,
	showing variation with~$t$ for given positions ($-5,-2,0,2,5)$.}
\end{figure}

As expected in this example, the limiting average position $\mathsf{E}X(t)$ fluctuates around~0, whereas the limiting variance $\mathsf{Var}X(t)$ is not (see Fig.~2). It
remains finite as the time becomes infinite.

\begin{figure}[!ht]
	\centering
	\includegraphics[scale=0.4]{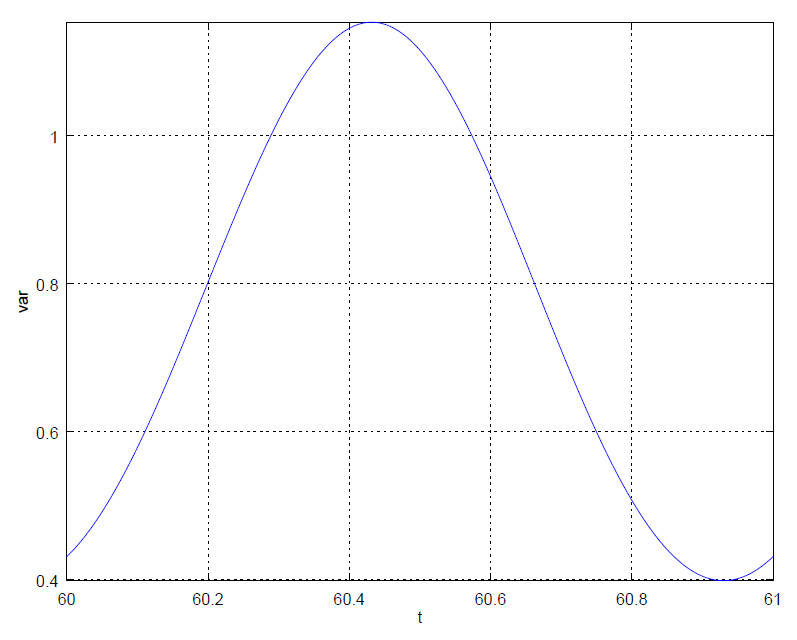}
	%\vspace{-1cm}
	\caption{Limiting variance $\mathsf{Var}X(t)$ of the particle position at time $t$.}
\end{figure}

%Consider a Markovian single-server queueing system
%with the queue of infinite capacity.

%\subsection*{\textbf{Example 2}}

As the second example we consider the double-ended queueing system
with the state space $\mathcal{X}= \{\cdots, -3, -2, -1, 0, 1, 2, 3, \cdots\} $.
Let $X(t)$ be the queue length of taxi or passenger at time $t$.
If $ X(t) > 0 $, the number of passengers in the system is $ X(t) $ and
there is no taxi queue. If $ X(t) < 0 $, the number of taxis in the system is $ -X(t) $ and there is no passenger
queue. If $ X(t) = 0 $, there is no taxi nor passenger.
Passengers and taxis arrive according to Poisson
process. Passengers (one to four passengers traveling together are considered as one passenger) arrive
to the queueing system according to a Poisson process with rate $ \la(t) $.
Obviously, $ \{X(t), t\ge 0\} $ is a one-dimensional continuous time
Markov chain.
The dynamic control of taxi depends on the state $X(t)$ of the system.
If there is no passenger (i.e. $ X(t) \le 0 $) waiting in the system,
the taxi arrival rate is $ \mu_1(t)$, otherwise (i.e. $ X(t) \ge 0 $) the taxi arrival
rate is $\mu_2(t)$ Obviously, the arrival rate of taxis with passengers is higher than that without passengers, i.e. $ \mu_1(t) \le \mu_2(t)$.
Passengers and taxis match according to the first-in-first-out discipline
and matching is instantaneous.
The transposed intensity matrix $A(t)$ for the considered problem 
has the following structure:

{\footnotesize
	\begin{equation*}
	A(t)=	\bordermatrix{
		~  &~ & -3 & -2 & -1 & 0 & 1 &2 &3 \cr
		~   & \cdots & \cdots & \cdots & \cdots & \cdots & \cdots & \cdots & \cdots & \cdots\cr
		-3 & \cdots& -\mu_{1}-\la &  \mu_{1} & 0 & 0 & 0 & 0 & 0 & \cdots \cr
		-2 & \cdots & \la &-\mu_{1}-\la & \mu_{1} & 0 & 0 & 0 & 0 & \cdots \cr
		-1& \cdots & 0& \la & -\mu_{1}-\la &  \mu_{1} & 0 & 0 & 0 &\cdots\cr
		0 & \cdots& 0 & 0 & \la & -\mu_{1}-\la & \mu_{2} & 0 & 0 &\cdots \cr
		1&  \cdots & 0& 0  & 0 & \la & -\mu_{2}-\la  &\mu_{2} & 0 &\cdots\cr
		2& \cdots& 0 &0 &0&0& \la & -\mu_{2}-\la &\mu_{2} &\cdots\cr
		3& \cdots & 0&0 &0&0& 0& \la & -\mu_{2}-\la &\cdots\cr
		~ & \cdots& \cdots & \cdots & \cdots & \cdots & \cdots & \cdots & \cdots & \cdots }.
	\end{equation*}
}

The \textit{Theorem~1} yields the upper bound
for the convergence rate, if a~double infinite sequence,
say $\{d_k, k= \pm 1, \pm 2 ,\dots\}$, can be found such that
$\int_0^\infty \beta^{**}(u)du=\infty$.
Let us put $ d_1=1, d_{-1} = c$ and $d_{k}=\delta^{k-1}, d_{-k}= c \cdot d_{k}$ for $k \ge 2$, where ${\delta < 1}$.
Then we have:
\begin{eqnarray}
\beta^{**}\left( t\right) =
\begin{cases} \left(\frac{1}{\delta}-1\right)\left(\delta\mu_1\left(t\right) - \la\left(t\right)\right)  , \ k < -1,\\
(1-\delta)\mu_{1}\left( t\right) + \left(1-\frac{1}{c}\right)\la\left(t\right), \ k = -1,\\
(1-\delta)\lambda\left( t\right) + \mu_{2}\left( t\right) - c\mu_1\left( t\right), \ k= 1,\\
\left(\frac{1}{\delta}-1\right)\left(\delta\la\left(t\right) - \mu_2\left(t\right)\right), \ k > 1.
\end{cases}
\label{2112}
\end{eqnarray}

\indent The value of ${\beta^{**}\left( t\right) =\min_{k \in \mathbb{Z} \setminus \{0\}}\beta^{**}_{k}\left( t\right)}$ cannot be written out
unless the exact values of $\lambda(t)$ and $\mu_k(t)$
are assumed. Let us fix
$ \lambda\left( t\right) = 2+\frac{\sin(2\pi t)}4 $, $ \mu_1\left( t\right)=1+\frac{\sin(2\pi t)}8 $ and $ \mu_2\left( t\right)=4+\frac{\cos(2\pi t)}4 $.
 Then if one puts
$ {d = \frac{8}{7}} $, ${ c = 2 }$
 and $d_{k}=\left(\frac{8}{7}\right)^{k-1}$ for ${k \ge 1}$, then the constants $\beta^{**}$ and $M$ from the \textit{Theorem~1} and \textit{Corollary~1} are equal to
$\beta^{**}=0.09375$ and $M=1$.
For the truncated process $X^*(t)$ with the truncation threshold $N=150$
one can put $ {d^* = \sqrt{2} }$,
${c^* = 2}$ and $d^*_{k}=\sqrt{2}^{k-1}$ for ${k \ge 1}$.
 Then the constants $\beta^{*}$
and $M^*$ from the~\textit{Theorem~2} are equal to $\beta^*=0.09375$ and $M^*=1$.
Thus, since  ${\|D(\vz(t)-\vz^*(t))\|\le 10^{-7}}$
from~\eqref{teor2} one gets
$$
\|\vp(t)-\vp^*(t)\| \le 10^{-7}.
$$
\noindent and from the comments, following the~\textit{Theorem~2},
one obtains
$$
\sum_{k \in \mathbb{Z} \setminus \{0\}}^\infty |k| |p_{k}(t)-p^*_{k}(t)|
\le 10^{-7}.
$$

Fig. 3 shows the variation of $p_{k}(t)$ with $t$ for
five different values $k$.

\begin{figure}[!ht]
	\centering
	\includegraphics[scale=0.4]{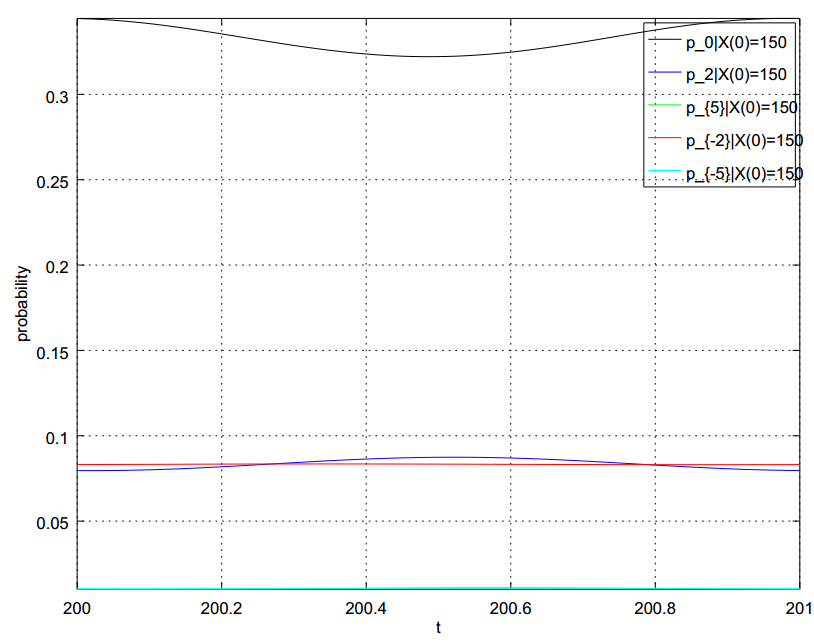}
	%\vspace{-1cm}
	\caption{Limiting probability of the process $X(t)$ at time $t$,
	showing variation with~$t$ for given positions ($-5,-2,0,2,5)$.}
\end{figure}

In this example, as in the previous one,
the limiting average position $\mathsf{E}X(t)$ fluctuates around~0
(see Fig.~4).

\begin{figure}[!ht]
	\centering
	\includegraphics[scale=0.4]{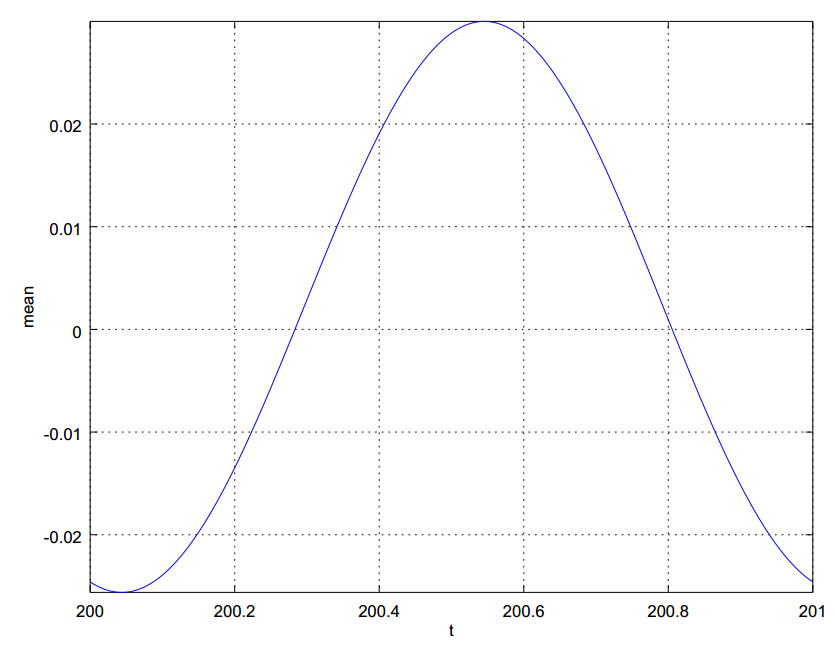}
	%\vspace{-1cm}
	\caption{Limiting expected value $\mathsf{E}X(t)$ of the process at time $t$.}
\end{figure}

The limiting variance $\mathsf{Var}X(t)$ is not around zero (see Fig.~5)
and remains finite as the time becomes infinite.

\begin{figure}[!ht]
	\centering
	\includegraphics[scale=0.4]{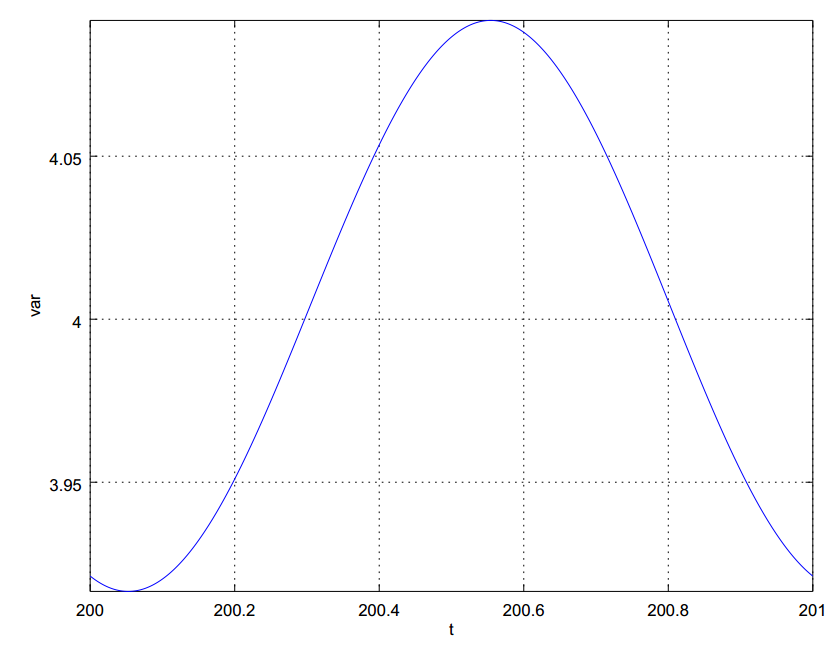}
	%\vspace{-1cm}
	\caption{Limiting variance $\mathsf{Var}X(t)$ of the process at time $t$.}
\end{figure}

%First of all this inequality is never negative.
%This inequality reminds of Markov's inequality, which can be used to bound
%the random variable $|X(t)|$:
%$$
%Pr(|X(t| >N) \le
%{E(|X(t)|) \over N},
%$$
%where $E(|X(t)|)$ denotes the average value of $|X(t)$ and is
%unknown. {\color {blue} ??? можно ли аналитически
%сравнить ${E(|X(t)|) \over N}$ и новую оценку???? }

%\bigskip

\section{\textbf{Conclusion}}

The developed theory for bilateral BDPs
facilitates their numerical analysis by providing upper
ergodicity and truncation bounds.
The latter can be used to understand when the limiting regime
is reached and show to properly truncate the
bi-infinite state space. The weak point of the
obtained results is the unknown bi-infinite sequence
of positive numbers $\{d_k\}$, for which no rule of thumb can be suggested 
and in each new use-case is has to be guessed.
Having no probabilistic meaning this sequence
can be considered as the analogue of Lyapunov~functions.

\section*{\textbf{Acknowledgements}}
This research was supported by Russian Science Foundation under grant 19-11-00020.

%\bigskip

\section*{\textbf{References}}
\renewcommand{\refname}{}

\end{document}